\documentclass{amsart}
\usepackage{graphicx}
\usepackage{latexsym}
\usepackage{amsfonts}
\usepackage[all]{xy}
\usepackage{amssymb, mathrsfs, amsfonts, amsmath}
\usepackage{amsbsy}
\usepackage{amsfonts}
\newtheorem*{acknowledgement}{Acknowledgement}
\newtheorem{corollary}{Corollary}

\newtheorem{definition}{Definition}

\newtheorem{proposition}{Proposition}
\newtheorem{remark}{Remark}
\newtheorem{theorem}{Theorem}

\numberwithin{equation}{section}

\begin{document}
\title[Anti-self dual quasi Yamabe Soliton]{A note on (anti-)self dual quasi Yamabe Soliton}
\author{Benedito Leandro Neto $^{1}$}
\address{$^{1}$ Departamento de Matem\'atica, Universidade de Bras\'ilia,
70910-900, Bras\'\i lia - DF, Brazil.}
 \email{B.L.Neto@mat.unb.br$^{1}$}
\thanks{$^{1}$ Supported by CNPq Proc. 149896/2012-3  Ministry of Science and Technology Brazil.}
%\date{September 16, 2014}

\begin{abstract}
In this note we prove that a (anti-)self dual quasi Yamabe soliton with positive sectional curvature is rotationally symmetric. This generalizes a recent result of G. Huang and H. Li in dimension four. Whence, (anti-) self dual gradient Yamabe solitons with positive sectional curvature is rotationally symmetric. We also prove that half conformally flat gradient Yamabe soliton has a special warped product structure provided that the potential function has no critical point.
\end{abstract}

\maketitle
\section{Introduction}

The Yamabe flow was introduced by R. Hamilton in an attempt to prove the Yamabe problem (see \cite{H}). Solitons are important to understand the geometric flow since they can appear as singularity models. Yamabe solitons are self-similar solutions for the Yamabe flow.

\begin{definition}
A Riemannian manifold $(M^n, g)$ with dimension $n\geq3$ is called a gradient Yamabe soliton if there exist a smooth potential function $f: M^n\rightarrow\mathbb{R}$ and a constant $\lambda$ such that
\begin{eqnarray}\label{yamabe}
(R-\lambda)g_{ij}=\nabla_{i}\nabla_{j}f.
\end{eqnarray}
\end{definition}
If $\lambda>0$, $\lambda=0$ or $\lambda<0$, then we have, respectively, a gradient Yamabe soliton shrinker, steady and expanding.

Motivated by the results on quasi Einstein manifolds (see \cite{C}), the theory of quasi Yamabe gradient solitons started to be investigated. Quasi Yamabe gradient solitons are generalizations of gradient Yamabe solitons.
\begin{definition}\label{definiçãoquasi}
A quasi Yamabe gradient soliton is a triple $(M^{n}, g, f)$, where $(M^{n}, g)$ is a Riemannian manifold of dimension $n\geq3$ with a smooth potential function  $f: M^n\rightarrow\mathbb{R}$ and two constants $\lambda, m$ ($m\neq0$) satisfying
\begin{equation}\label{quasi}
(R-\lambda)g_{ij}=\nabla_{i}\nabla_{j}f-\frac{1}{m}\nabla_{i}f\nabla_{j}f.
\end{equation}
\end{definition} When $m\rightarrow\infty$ then (\ref{quasi}) reduces to (\ref{yamabe}) and if $f$ is constant we say that the quasi Yamabe gradient soliton is trivial.

Recently $4$-dimensional manifolds has been widely studied (see \cite{BLR}, \cite{br}, \cite{C}, \cite{CW} and \cite{HL}). In what follows $M^4$ will denote an oriented 4-dimensional manifold and $g$ is a Riemannian metric on $M^4.$ We emphasize that 4-manifolds are fairly special.  For instance, following the notations used in \cite{handbook} (see also \cite{ku} and \cite{besse}), given any local orthogonal frame $\{e_{1}, e_{2}, e_{3}, e_{4}\}$ on open set of $M^4$ with dual basis $\{e^{1}, e^{2}, e^{3}, e^{4}\},$ there exists a unique bundle morphism $\ast$ called {\it Hodge star} (acting on bivectors), such that $$\ast(e^{1}\wedge e^{2})=e^{3}\wedge e^{4}.$$ This implies that $\ast$ is an involution, i.e. $\ast^{2}=Id.$ In particular, this implies that the bundle of $2$-forms on a 4-dimensional oriented Riemannian manifold can be invariantly decomposed as a direct sum $\Lambda^2=\Lambda^{2}_{+}\oplus\Lambda^{2}_{-}.$ This allows us to conclude that the Weyl tensor $W$ is an endomorphism of $\Lambda^2=\Lambda^{+} \oplus \Lambda^{-} $ such that
\begin{equation}
\label{decW}
W = W^+\oplus W^-.
\end{equation} We recall that $dim_{\Bbb{R}}(\Lambda^2)=6$ and $dim_{\Bbb{R}}(\Lambda^{\pm})=3.$ Also, it is well-known that
\begin{equation}
\label{6A}
\Lambda^{+}=span\Big\{\frac{e^{1}\wedge e^{2}+e^{3}\wedge e^{4}}{\sqrt{2}},\,\frac{e^{1}\wedge e^{3}+e^{4}\wedge e^{2}}{\sqrt{2}},\,\frac{e^{3}\wedge e^{2}+e^{4}\wedge e^{1}}{\sqrt{2}}\Big\}
\end{equation}
 and
\begin{equation}
\label{6B}
 \Lambda^{-}=span\Big\{\frac{e^{1}\wedge e^{2}-e^{3}\wedge e^{4}}{\sqrt{2}},\,\frac{e^{1}\wedge e^{3}-e^{4}\wedge e^{2}}{\sqrt{2}},\,\frac{e^{3}\wedge e^{2}-e^{4}\wedge e^{1}}{\sqrt{2}}\Big\}.
\end{equation} From this, the bundles $\Lambda^{+}$ and $\Lambda^{-}$ carry natural orientations such that the bases (\ref{6A}) and (\ref{6B}) are positive-oriented.

The decomposition of the Weyl tensor on $4$-dimensional manifolds allow us to deduce the following equation
\begin{equation}
\label{W+}
W_{\,p\,q\,r\,s}^{+}=\frac{1}{2}\big(W_{p\,q\,r\,s}+W_{\overline{p}\,\overline{q}\,r\,s}\big),
\end{equation} where $(\overline{p}\,\overline{q}),$ for instance, stands for the dual of $(p\,q),$ that is, $(\overline{p}\,\overline{q}\,p\,q)=\sigma(1234)$ for some even permutation $\sigma$ in the set $\{1,2,3,4\}$ (cf. Equation 6.17, p. 466 in \cite{handbook}). In particular, we have $$W_{1234}^{+}=\frac{1}{2}\big(W_{1234}+W_{1212}\big).$$ For more details we refer to \cite{handbook}, \cite{ku} and \cite{besse}. We say that a $4$-dimensional manifold is half conformally flat if it is self or anti-self dual, namely if $W^{-}=0$ or $W^{+}=0$.

Inspired by \cite{CW} and \cite{HL}, we consider $4$-dimensional quasi Yamabe gradient solitons. More precisely, we prove:
\begin{theorem}\label{teoremaprincipal}
Let $(M^4, g, f)$ be a nontrivial complete anti-self dual (or self dual) quasi Yamabe gradient soliton satisfying (\ref{quasi})
with positive sectional curvature. Then $(M^4, g, f)$ is rotationally symmetric.
\end{theorem}
We have the following result from Theorem \ref{teoremaprincipal} immediately:
\begin{corollary}\label{coro1}
Let $(M^4, g, f)$ be a nontrivial complete anti-self dual (or self dual) gradient Yamabe soliton satisfying (\ref{yamabe})
with positive sectional curvature. Then $(M^4, g, f)$ is rotationally symmetric.
\end{corollary}

In fact, we can improve corollary \ref{coro1}.

\begin{theorem}\label{teoremaprincipal1}
Let $(M^4, g, f)$ be a nontrivial complete anti-self dual (or self dual) gradient Yamabe soliton satisfying (\ref{yamabe}) and suppose that $f$ has no critical points. Then $(M^4, g, f)$ is the warped product
\begin{eqnarray*}
(\mathbb{R}, dr^{2})\times_{|\nabla f|}(N^{3}, \bar{g}),
\end{eqnarray*}
where $(N^{3}, \bar{g})$ is a space form (i.e, of constant sectional curvature).
\end{theorem}

\begin{remark}
In \cite{DS}, Daskalopoulos and Sesum proved that gradient Yamabe solitons with positive sectional curvature are rotationally symmetric under the assumption that the metric $g$ is locally conformally flat. Cao, Sun and Zhang \cite{CSZ} proved that gradient Yamabe solitons has a special warped product structure.
\end{remark}

\section{Preliminares}

The Weyl tensor $W$ is defined by the following decomposition formula
\begin{eqnarray}\label{weyl}
R_{ijkl}&=&W_{ijkl}+\frac{1}{n-2}\big(R_{ik}g_{jl}+R_{jl}g_{ik}-R_{il}g_{jk}-R_{jk}g_{il}\big)\nonumber\\
&-&\frac{R}{(n-1)(n-2)}\big(g_{jl}g_{ik}-g_{il}g_{jk}\big).
\end{eqnarray}
It is well know that $W=0$ for $n=3$. A $4$-dimensional manifold is locally conformally flat if and only if $W=0$.

In \cite{CC}, Cao and Chen defined the tensor $D$. This tensor is the link between the Weyl tensor and quasi Yamabe solitons (see proposition \ref{prop1}). We define the $3$-tensor $D$ by
\begin{eqnarray}\label{tensorD}
D_{ijk}&=&\frac{1}{n-2}(R_{jk}\nabla_{i}f-R_{ik}\nabla_{j}f)+\frac{1}{(n-1)(n-2)}(R_{il}\nabla^{l}fg_{jk}-R_{jl}\nabla^{l}fg_{ik})\nonumber\\
&-&\frac{R}{(n-1)(n-2)}(g_{jk}\nabla_{i}f-g_{ik}\nabla_{j}f).
\end{eqnarray}
The tensor $D$ is skew-symmetric in their first two indices and trace-free, i.e,
\begin{eqnarray*}
D_{ijk}=-D_{jik}\quad\mbox{and}\quad g^{ij}D_{ijk}=g^{ik}D_{ijk}=g^{jk}D_{ijk}=0.
\end{eqnarray*}
\begin{proposition}\cite{HL}\label{prop1}
Let $(M^{n}, g, f)$ be a nontrivial complete quasi Yamabe gradient soliton satisfying (\ref{quasi}). Then the $3$-tensor $D$ is related to the Weyl curvature tensor by
\begin{eqnarray*}
D_{ijk}=W_{ijkl}\nabla^{l}f.
\end{eqnarray*}
\end{proposition}

\begin{theorem}\cite{HL}\label{theorem1}
Let $(M^{4}, g, f)$ be a nontrivial complete quasi Yamabe gradient soliton satisfying (\ref{quasi}) with positive sectional curvature and $D_{ijk}=0$. Then $(M^{4}, g, f)$ is rotationally symmetric.
\end{theorem}

\begin{theorem}\cite{CSZ}\label{theoremCSZ1}
Let $(M^{n}, g, f)$ be a nontrivial complete gradient Yamabe soliton satisfying equation (\ref{yamabe}). Then $|\nabla f|^{2}$ is constant on regular level surfaces of $f$, and either

(i) $f$ has a unique critical point at some point $x_{0}\in M^{n}$, and $(M^{n}, g, f)$ is rotationally symmetric and equal to  the warped product
\begin{eqnarray*}
([0,+\infty), dr^{2})\times_{|\nabla f|}(\mathbb{S}^{n-1}, \bar{g}_{can}),
\end{eqnarray*}
where $\bar{g}_{can}$ is the metric on $\mathbb{S}^{n-1}$, or

(ii) $f$ has no critical point and $(M^{n}, g, f)$ is the warped product
\begin{eqnarray}\label{structure}
(\mathbb{R}, dr^{2})\times_{|\nabla f|}(N^{n-1}, \bar{g}),
\end{eqnarray}
where $(N^{n-1}, \bar{g})$ is a Riemannian manifold of constant scalar curvature.
\end{theorem}

When $(M^{n}, g, f)$ is locally conformally flat we have

\begin{theorem}\cite{CSZ}\label{theoremCSZ}
Let $(M^{n}, g, f)$ be a nontrivial complete gradient Yamabe soliton satisfying equation (\ref{yamabe}). Suppose $f$ has no critical point and is locally conformally flat. Then $(M^n, g, f)$ is the warped product
\begin{eqnarray*}
(\mathbb{R}, dr^{2})\times_{|\nabla f|}(N^{n-1}, \bar{g}),
\end{eqnarray*}
where $(N^{n-1}, \bar{g})$ is a space form (i.e, of constant sectional curvature).
\end{theorem}

In this paper we show that the hypothesis locally conformally flat in theorem \ref{theoremCSZ} can be replaced by the weaker condition (anti-)self dual Weyl tensor (or half conformally flat) on four dimensional case.

\begin{remark}
It is worth point out that compact quasi Yamabe gradient solitons and compact Yamabe gradient solitons has constant scalar curvature (see \cite{DS} and \cite{HL}).
\end{remark}

\section{Proof of the main results}

\subsection{Proof of Theorem \ref{teoremaprincipal}}
\begin{proof}
Now, we follow the same ideas from \cite{CW}. Consider a $4$-dimensional Riemannian manifold satisfying definition \ref{definiçãoquasi} with
 $W^{+}=0$, i.e, a complete anti-self dual quasi Yamabe gradient soliton. Then, from (\ref{W+}) we have
\begin{eqnarray*}
W_{i\,j\,k\,l}+W_{\overline{i}\,\overline{j}\,k\,l}=0.
\end{eqnarray*}
Whence,
\begin{eqnarray*}
W_{i\,j\,k\,l}\nabla^{l}f+W_{\overline{i}\,\overline{j}\,k\,l}\nabla^{l}f=0.
\end{eqnarray*}
This implies from proposition \ref{prop1} that
\begin{eqnarray*}
D_{i\,j\,k}+D_{\overline{i}\,\overline{j}\,k}=0.
\end{eqnarray*}

Now, for some oriented orthonormal basis $\{e_{i}\}_{i=1}^{4}$ diagonalizing the Ricci tensor with associated eigenvalues $\mu_{k}$, $k=1,\ldots,4$, respectively, i.e, $R_{ij}=\mu_{i}\delta_{ij}$, we have
\begin{eqnarray}\label{quadrocomutativo}
D_{12k}+D_{34k}=0,\quad D_{13k}+D_{42k}=0\quad\mbox{and}\quad D_{14k}+D_{23k}=0.
\end{eqnarray}
From (\ref{tensorD}) we get
\begin{eqnarray}\label{tensorD1}
D_{ijk}=0,\quad{for}\quad i\neq j\neq k.
\end{eqnarray}
Therefore, from (\ref{quadrocomutativo}) and (\ref{tensorD1}) we obtain
\begin{eqnarray}\label{tensorD2}
D_{iji}=0,\quad\mbox{for}\quad i,j=1,\ldots,4.
\end{eqnarray}
We also have that the tensor $D$ is skew-symmetric, i.e, $D_{iij}=0$. Finally we get that $D\equiv0$. Since the sectional curvature is positive we can apply theorem \ref{theorem1} to get the result.
\end{proof}

\subsection{Proof of Theorem \ref{teoremaprincipal1}}
\begin{proof}
Let $(M^4, g, f)$ be a complete anti-self dual (or half conformally flat) gradient Yamabe soliton such that $f$ has no critical point. From theorem \ref{theoremCSZ1}-(ii) we get that $(M^4, g, f)$ has a warped product structure (\ref{structure}).
We want to prove that $(N^{3}, \bar{g})$ is a space form. For this, it suffices to show that $\bar{g}$ is Einstein. Since three dimensional Einstein manifolds has constant sectional curvature. To conclude our result we need to know the warped product structure of complete gradient Yamabe soliton (see \cite{CSZ}).

By theorem \ref{teoremaprincipal} and proposition \ref{prop1} we already know that $W_{ijkl}\nabla^{l}f=0$, i.e,
\begin{eqnarray}\label{blz}
W(\cdot,\cdot,\cdot,\nabla f)=0.
\end{eqnarray}
Consider the level surface $\Sigma=f^{-1}(c)$ where $c$ is any regular value of the potential function $f$. Suppose that $I$ is an open interval containing $c$ such that $f$ has no critical point in the open neighborhood $U_{I}=f^{-1}(I)$. Since $|\nabla f|^{2}$ is constant on $\Sigma$, we can make a chance of variable so that we can express the metric $g$ in $U_{I}$ as in theorem \ref{theoremCSZ1}-(ii) (see \cite{CSZ}) .
Fix any local coordinates system $\theta=(\theta_{2},\theta_{3},\theta_{4})$ on $N^{3}$, and choose $(x_{1}, x_{2}, x_{3}, x_{4})=(r, \theta_{2},\theta_{3},\theta_{4})$. Let $\nabla r=\frac{\partial}{\partial r}$, then $|\nabla r|=1$ and $\nabla f=f'(r)\frac{\partial}{\partial r}$ on $U_{I}$. Whence, we get from (\ref{blz})
\begin{eqnarray*}
0=W(\cdot,\cdot,\cdot,\nabla f)=f'(r)W(\cdot,\cdot,\cdot,\nabla r).
\end{eqnarray*}
Since $f$ has no critical point we get
\begin{eqnarray}\label{blz1}
W(\cdot,\cdot,\cdot,\nabla r)=0.
\end{eqnarray}
Therefore, from arguments on \cite{CSZ} (see theorem 1.4), the Weyl tensor formula for an arbitrary warped product dimensional manifold give us
\begin{eqnarray}\label{Weylwarped}
2W(\nabla r,\theta_{a},\nabla r,\theta_{b})=\frac{\bar{R}}{3}\bar{g}(\theta_{a},\theta_{b})-\bar{R}ic(\theta_{a},\theta_{b}),\quad\mbox{for}\quad a, b\in\{2, 3, 4\},
\end{eqnarray}
where $\bar{R}$ and $\bar{R}ic$ stand, respectively, for the scalar curvature and Ricci tensor for $(N^{3}, \bar{g})$. Therefore, from (\ref{blz1}) and (\ref{Weylwarped}) we have theorem \ref{teoremaprincipal1}.
\end{proof}

\begin{acknowledgement}
The author was partially supported by grant from CNPq-Brazil.
%Moreover, he wants to thank Professor E. Ribeiro Jr. for valuable discussion about this issue.
%Finally, the author wants to thank the referee for his careful reading and helpful suggestions.
\end{acknowledgement}

\end{document}